\documentclass[10pt,twoside]{siamltex}
\usepackage{amsmath}
\usepackage{amsfonts}
\newcommand{\R}{{\mathbb R}}
\newcommand{\N}{{\mathbb N}}

\begin{document}
\title{Positive eigenvalues and two-letter generalized words
\thanks{Received by the editors on .
Accepted for publication on . Handling Editor: .}}
\author{C. Hillar\thanks{\mbox{Mathematics Department,
University of California, Berkeley Berkeley, CA 94720, USA}
(chillar@uclink.berkeley.edu). Supported under a National Science
Foundation Graduate Research Fellowship.} \and C. R.
Johnson\thanks{\mbox{Department of Mathematics, College of William
\& Mary, Williamsburg, VA 23187-8795, USA} (crjohnso@math.wm.edu,
ilya@math.wm.edu). This research was supported by NSF REU Grant
DMS 99-87803.} \and I. M. Spitkovsky \footnotemark[3]}

\markboth{C. Hillar, C. Johnson and I. Spitkovsky} {Two letter
words}

\maketitle

\begin{abstract}
A generalized word in two letters $A$ and $B$ is an expression of
the form $W=A^{\alpha_1}B^{\beta_1}A^{\alpha_2}B^{\beta_2}\cdots
A^{\alpha_N}B^{\beta_N}$ in which the exponents $\alpha_i$,
$\beta_i$ are nonzero real numbers. When independent positive
definite matrices are substituted for $A$ and $B$, we are
interested in whether $W$ necessarily has positive eigenvalues.
This is known to be the case when $N=1$ and has been studied in
case all exponents are positive by two of the authors. When the
exponent signs are mixed, however, the situation is quite
different (even for 2-by-2 matrices), and this is the focus of the
present work.
\end{abstract}

\begin{keywords}
Positive definite matrices, Projections, Generalized word.
\end{keywords}
\begin{AMS}
15A18, 15A57 \end{AMS}

Let $A,B$ be positive definite $n\times n$ matrices. Then, as is
well known \cite[p. 465]{HJ1}, the eigenvalues of the product $AB$
are real and positive. Moreover, for all $\alpha,\beta\in\R$ the
matrices $A^\alpha$ and $B^\beta$ are positive definite together
with $A,B$. Thus, the eigenvalues of $A^\alpha B^\beta$ are real
and positive as well.

In this paper, we are concerned with possible generalizations of
this simple observation to products
$W(A,B)=A^{\alpha_1}B^{\beta_1}A^{\alpha_2}\ldots$ Such
expressions, when the $\alpha$'s and $\beta$'s are positive
integers, have been studied in \cite{HiJo1} and when  $\alpha$'s
and $\beta$'s are positive reals in subsequent work. Applying an
appropriate similarity if necessary, we may without loss of
generality suppose that $W(A,B)$ ends with a power of $B$. In
other words,
\begin{equation}\label{word}
W(A,B)= A^{\alpha_1}B^{\beta_1}A^{\alpha_2}B^{\beta_2}\cdots
A^{\alpha_N}B^{\beta_N}\end{equation}
($\alpha_j,\beta_j\in\R\setminus\{ 0\})$. We will say that
(\ref{word}) is a {\em generalized word} (g-word) in $A,B$ of
class $N$.

{\bf Problem}. {\em Under what additional conditions on $A,B$
and/or the structure of the g-word {\em (\ref{word})} is it true
that all the eigenvalues of $W(A,B)$ are positive?}

The above observation means that there are no additional
conditions on $A$ and $B$ for $g$-words of class 1. Another
trivial sufficient condition is the commutativity of $A$ and $B$
(which holds, in particular, for $n=1$). 
Starting with $n=N=2$, it is easy to give examples of g-words
(\ref{word}) with positive definite $A$, $B$ and the spectrum not
lying in $\R_+$. The simplest such word is $ABA^{-1}B^{-1}$. That
this word does not guarantee positive spectrum can be seen from
the following, more precise, statement.
\begin{theorem}\label{th:four}
Let $A$ have exactly two distinct eigenvalues. Then the spectrum
of $A^mBA^{-m}B^{-1}$ is positive for all $m\in\N$ if and only if
$A$ and $B$ commute.
\end{theorem}
\begin{proof}
Using a unitary similarity if necessary, we may put $A$ in the
form
\[
A=\left[\begin{array}{cc} \lambda_1 I_{n_1} & 0 \\ 0 &\lambda_2
I_{n_2}\end{array}\right],\] where $\lambda_1>\lambda_2>0$; denote
the respective partition of $B$ by
\[
B=\left[\begin{array}{cc} B_{11} & B_{12} \\ B_{21} &
B_{22}\end{array}\right]\] (due to self adjointness of $B$, the
blocks $B_{11}$, $B_{22}$ also are self adjoint, and
$B_{21}=B_{12}^*$). Then
\[
A^mBA^{-m}=\gamma^{-m}\left[\begin{array}{cc}\gamma^mB_{11} & B_{12}\\
\gamma^{2m} B_{21} & \gamma^m B_{22}\end{array}\right],
\]
where $\gamma=\lambda_2/\lambda_1<1$. Thus, there exists the limit
of $\gamma^mA^mBA^{-m}B^{-1}$ when $m\to\infty$, and this limit
equals
\begin{equation}\label{limit}
\left[\begin{array}{cc} -B_{12}C^{-1}B_{21}B_{11}^{-1} & B_{12}C^{-1}\\
0 & 0\end{array}\right],
\end{equation}
where $C=B_{22}-B_{21}B_{11}^{-1}B_{12}$ is positive definite due
to the positive definiteness of $B$ (see, e.g., \cite[p.
475]{HJ1}). Suppose that the eigenvalues of all the matrices
$A^mBA^{-m}B^{-1}$ are positive. Then all the eigenvalues of the
left upper block of the matrix (\ref{limit}) are non-negative. In
other words, the spectrum of the matrix
$B_{12}C^{-1}B_{21}B_{11}^{-1}$ must be non-positive. The latter
being a product of a non-negative definite matrix
$B_{12}C^{-1}B_{21}$ and a positive definite matrix $B_{11}^{-1}$,
this is only possible if it is the zero matrix. But then
$0=B_{12}C^{-1}B_{21}=(B_{12}C^{-1/2})(B_{12}C^{-1/2})^*$, so that
$B_{12}=0$. This implies that $B_{21}=0$ as well. In other words,
$B$ commutes with $A$.
\end{proof}

Observe that in Theorem~\ref{th:four} both $A$ and $B$ appear with
powers of different sign, and that for the g-words of class 1 this
situation is impossible. So, it is natural to entertain a
conjecture that g-words (\ref{word}) with powers of the same sign
have positive spectra. As it happens, this is also not true (even
for words of class 2 and natural exponents) but the respective
example is much harder to come by. The simplest known example of
this kind is the word $ABA^2B^2$, with
\begin{equation}\label{aba2b2}
A=\left[\begin{array}{ccc} 1& 20& 210\\
20 & 402 & 4240 \\ 210 & 4240 & 44903 \end{array}\right],
\quad B= \left[\begin{array}{ccc} 36501 & -3820 & 190 \\
-3820 & 401 & -20 \\ 190 & -20 & 1
\end{array}\right]\end{equation}
(see \cite{HiJo1}). Note that in (\ref{aba2b2}) $n=3$ and all the
eigenvalues of both matrices $A$ and $B$ are distinct. The next
two theorems show that these features are indeed necessary for
such an example. Let us prove an auxiliary statement first.
\begin{lemma}\label{l:multeig}
Let one of the matrices $A$, $B$ have an eigenvalue of
multiplicity at least $n-1$ and in {\em (\ref{word})} all the
powers of the other matrix be of the same sign. Then $W(A,B)$ has
at least one positive eigenvalue.
\end{lemma}
\begin{proof}
Without loss of generality (by a simple change of notation if
necessary) we may suppose that $A$ is the matrix with an
eigenvalue $\lambda_1$ of multiplicity $n-1$; denote its remaining
eigenvalue by $\lambda_2$. Switching from $B$ to $B^{-1}$ if
necessary, we may also suppose that $\beta_1,\ldots,\beta_N\geq
0$.

{\sl Case 1}. $\beta_1,\ldots,\beta_N$ are integers. Let $U$ be a
unitary similarity diagonalizing $A$:
\begin{equation}\label{A0}
A_0=U^*AU=\left[\begin{array}{ccccc}\lambda_1 & 0 & \ldots & 0 & 0 \\
0 & \lambda_1 &  & 0 & 0\\ \vdots & &\ddots & \vdots & \vdots \\
0 & \ldots & 0 & \lambda_1 & 0\\0 & 0 & \ldots & 0 &
\lambda_2\end{array}\right] .
\end{equation}
By an appropriate choice of $U$ (which consists in multiplying the
original one on the right by $V\oplus [1]$, where $V$ is some
$(n-1)\times (n-1)$ unitary matrix), we may suppose that the left
upper $(n-1)\times (n-1)$ block of $B$ also is diagonalized.
Multiplying $V$ on the right by a diagonal unitary matrix with
suitably chosen arguments of its diagonal entries, we can force
all the elements of the last column in $B_0=U^*BU$ to become
non-negative. But then all elements of its last row automatically
become non-negative as well. In other words, simultaneously with
(\ref{A0}) the following decomposition also holds:
\begin{equation}\label{B0}
B_0=U^*BU=\left[\begin{array}{ccccc}\mu_1 & 0 & \ldots & 0 & \gamma_1 \\
0 & \mu_2 &  & 0 & \gamma_2 \\ \vdots & &\ddots & \vdots & \vdots \\
0 & \ldots & 0 & \mu_{n-1} & \gamma_{n-1}\\ {\gamma_1} &
{\gamma_2} & \ldots & {\gamma_{n-1}} & \ \mu_n
\end{array}\right] .
\end{equation}
Both matrices $A_0$ and $B_0$ are (entry-wise) non-negative. 
Thus, $W(A_0,B_0)$ also is entry-wise non-negative, and (at least)
one of its eigenvalues is positive due to Perron's theorem. 
But $W(A_0,B_0)=U^*W(A,B)U$, and the result follows.

{\sl Case 2}. $\beta_1,\ldots,\beta_N$ are rational. Let $Q (\in
\N)$ be their least common denominator. Considering $B^{1/Q}$, we
reduce this situation to Case 1.

{\sl Case 3}. Arbitrary (non-negative) $\beta_1,\ldots,\beta_N$.
For each $j=1,\ldots,N$, introduce a sequence $\beta_j^{(k)}$ of
non-negative rational numbers such that
$\lim_{k\to\infty}\beta_j^{(k)}=\beta_j$. Let
\[
W_k(A,B)=A^{\alpha_1}B^{\beta_1^{(k)}}A^{\alpha_2}B^{\beta_2^{(k)}}\cdots
A^{\alpha_N}B^{\beta_N^{(k)}}.
\]
Then each of the matrices $W_k(A,B)$ has a positive eigenvalue
(due to Case 2), and their limit $W(A,B)$ is invertible. From
continuity considerations it follows that $W(A,B)$ also has a
positive eigenvalue.
\end{proof}
\begin{theorem}\label{th:n=2} Let $n=2$, and let all powers of either
$A$ or $B$ in {\em (\ref{word})} be of the same sign. Then all the
eigenvalues of $W(A,B)$ are positive.
\end{theorem}
\begin{proof}
Since $n-1=1$, both $A$ and $B$ have eigenvalues of multiplicity
$n-1$. Hence, conditions of Lemma~\ref{l:multeig} are satisfied,
so that at least one eigenvalue of $W(A,B)$ is positive. But the
product of the two eigenvalues, $\det W(A,B)$, is positive as
well. Thus, the second eigenvalue is also positive.
\end{proof}
\begin{theorem}\label{th:n=3}
Let $n=3$, and suppose that at least one of the matrices $A$, $B$
has a multiple eigenvalue. If all the powers of the other matrix
in {\em (\ref{word})} are of the same sign, then all the
eigenvalues of $W(A,B)$ are positive.
\end{theorem}
\begin{proof}
Since $n-1=2$, conditions of Lemma~\ref{l:multeig} are met. We
will use representations (\ref{A0}), (\ref{B0}) from its proof,
which in case $n=3$ take the form
\[
A_0=\left[\begin{array}{ccc}\lambda_1 & 0 & 0\\ 0 & \lambda_1 & 0 \\
0 & 0 & \lambda_2\end{array}\right], \quad
B_0=\left[\begin{array}{ccc}\mu_1 & 0 & \gamma_1 \\ 0 & \mu_2 & \gamma_2 \\
\gamma_1 & \gamma_2 & \mu_3\end{array}\right].
\]
If $\gamma_1=0$ or $\gamma_2=0$ then $A_0$ and $B_0$ are
simultaneously in the block diagonal form, so that $W(A_0,B_0)$ is
a direct sum of a positive scalar and $W(A_1,B_1)$, where $A_1$
and $B_1$ are $2\times 2$ positive definite matrices. The result
then follows from Theorem~\ref{th:n=2}.

If both $\gamma_1$ and $\gamma_2$ are strictly positive, we will
again consider first the case of natural powers of $B$. There is
no need to consider the case $N=1$; in all other cases
$W(A_0,B_0)$ is entry-wise positive. According to Perron's
theorem, its positive eigenvalue $\eta_1$ coinciding with the
spectral radius is the only eigenvalue of this magnitude. Thus,
$\eta_1$ is the eigenvalue of $W(A,B)$ and the other two
eigenvalues satisfy $|\eta_3|\leq |\eta_2|<\eta_1$. Observe now
that $W(A,B)^{-1}$ is a word in $A^{-1}$, $B^{-1}$, and that
$A^{-1}$, $B^{-1}$ satisfy conditions of Lemma~\ref{l:multeig}
simultaneously with $A$, $B$. Thus, the biggest by its absolute
value eigenvalue $\eta_3^{-1}$ of $W(A,B)^{-1}$ must be positive
as well. From this, and the positivity of $\det
W(A,B)=\eta_1\eta_2\eta_3$ we conclude that the remaining
eigenvalue $\eta_2$ is also positive.

The case of arbitrary real $\beta_1,\ldots,\beta_N$ of the same
sign can be now covered in exactly the same manner as in the proof
of Lemma~\ref{l:multeig}.
\end{proof}

Our next result shows that in Theorem~\ref{th:n=2} it is not the
size of the matrices that counts but actually the number of their
distinct eigenvalues.
\begin{theorem}\label{th:2eig}
Suppose that each of the matrices $A$ and $B$ has at most two
distinct eigenvalues and that in {\em (\ref{word})} all powers of
either $A$ or $B$ are of the same sign. Then, for an arbitrary
$n$, all the eigenvalues of $W(A,B)$ are positive.
\end{theorem}
\begin{proof}
If $\lambda_1$ and $\lambda_2$ are the only eigenvalues of $A$,
then $A=(\lambda_1-\lambda_2)P+\lambda_2 I$, where $P$ is a
certain orthoprojection. Similarly, $B=(\mu_1-\mu_2)Q+\mu_2 I$,
where $Q$ is another orthoprojection. It is well known (see, e.g.,
\cite{Dav}, \cite{Dix}, or \cite{Hal69}) that, for any two
orthoprojections $P$ and $Q$, there is a unitary similarity $U$
such that
\begin{equation}\label{proj} P_0= U^*PU=P_1\oplus
P_2\oplus\cdots\oplus P_N, \ Q_0= U^*QU=Q_1\oplus
Q_2\oplus\cdots\oplus Q_N,
\end{equation}
where the size of $P_j$ is the same as the size of $Q_j$ and does
not exceed 2 ($j=1,\ldots,N$). But then
\[
U^*W(A,B)U=W(A_1,B_1)\oplus W(A_2,B_2)\oplus\cdots\oplus
W(A_N,B_N),
\]
where $A_j=(\lambda_1-\lambda_2)P_j+\lambda_2 I$,
$B_j=(\mu_1-\mu_2)Q_j+\mu_2 I$ are either positive numbers or
positive definite $2\times 2$ matrices. Due to
Theorem~\ref{th:n=2}, the eigenvalues of $W(A_j,B_j)$ are all
positive. The same is true for their direct sum $U^*W(A,B)U$, and
thus for $W(A,B)$ itself.
\end{proof}

Let us say that the sequence
$\alpha_1,\beta_1,\ldots,\alpha_N,\beta_N\, (\in (\R\setminus\{
0\})^{2N}$) is {\em 2-good} if the word (\ref{word}) has positive
eigenvalues for all positive definite $2\times 2$ matrices $A,B$.
Of course, {\em $k$-good} sequences can be defined in a similar
way for any $k\in\N$, and every $k$-good sequence is also $j$-good
for $j<k$.  According to Theorem~\ref{th:2eig}, any sequence for
which either all $\alpha$'s or all $\beta$'s are of the same sign
is 2-good. Many such sequences are $k$-good for all positive
integers $k$, as discussed in \cite{HiJo1}. On the other hand,
Theorem~\ref{th:four} implies that the sequence
$\alpha,\beta,-\alpha,-\beta$ is not 2-good. In fact, the
magnitudes of the exponents are in this case irrelevant: any
sequence $\alpha_1,\beta_1,\alpha_2,\beta_2$ with
$\alpha_1\alpha_2<0$, $\beta_1\beta_2<0$ is not 2-good. This
statement is a particular case of a more general one, the
formulation of which requires some preparation.

Consider the following {\em cancellation rule} for the sequences
$\alpha_1,\beta_1,\ldots,\alpha_m,\beta_m$, $m\in\N$: if
$\alpha_j\alpha_{j+1}>0$ for some $j\in\{ 1,\ldots,m\}$ (where by
convention $\alpha_{m+1}=\alpha_1$), then $\alpha_j, \beta_j$ are
omitted from the sequence. Similarly, if $\beta_j\beta_{j+1}>0$
then $\alpha_{j+1}, \beta_{j+1}$ are omitted. The sequence
$\alpha_1,\beta_1,\ldots,\alpha_m,\beta_m$ is {\em irreducible} if
no cancellations (in the above sense) are possible. Observe that
the signs of both $\alpha_1,\alpha_2,\ldots$  and
$\beta_1,\beta_2,\ldots$ in an irreducible sequence alternate. We
will say that $m$ is the {\em reduced} class of the sequence
$\alpha_1,\beta_1,\ldots,\alpha_N,\beta_N$ if there is an
irreducible sequence consisting of $2m$ terms obtained from
$\alpha_1,\beta_1,\ldots,\alpha_N,\beta_N$ by a repeated
application of the cancellation rule.

\begin{theorem}\label{th:not2good}
Any sequence $\alpha_1,\beta_1,\ldots \alpha_N,\beta_N$ of the
reduced class $m \equiv 2$ or $3$ {\em mod} $4$ is not 2-good.
\end{theorem}
\begin{proof}
Switching from $A$ to $A^{-1}$ and/or from $B$ to $B^{-1}$ if
necessary, we may without loss of generality suppose that the
first $\alpha$ and $\beta$ remaining after the cancellation
procedure are both positive. Then let
\[ A=\left[\begin{array}{cc}1 & 0\\ 0 & \epsilon\end{array}\right],
\quad B=\left[\begin{array}{cc}1/2+\epsilon & 1/2 \\ 1/2 & 1/2
\end{array}\right] \]
for some $\epsilon>0$. An easy computation shows that the matrix
\begin{equation}\label{epsilon}
2^{\sum_{\beta_j<0}\beta_j}\epsilon^{-(\sum_{\alpha_j<0}\alpha_j+
\sum_{\beta_j<0}\beta_j)}A^{\alpha_1}B^{\beta_1}A^{\alpha_2}B^{\beta_2}\cdots
A^{\alpha_N}B^{\beta_N}\end{equation} is the product of $2N$
matrices the $(2j-1)$-st of which is $\left[\begin{array}{cc}1 & 0\\
0 & \epsilon^{\alpha_j}\end{array}\right]$ if $\alpha_j>0$ and
$\left[\begin{array}{cc}\epsilon^{-\alpha_j} & 0\\ 0 &
1\end{array}\right]$ if $\alpha_j<0$, and the $2j$-th of which is
$\left[\begin{array}{cc}1/2+\epsilon & 1/2 \\ 1/2 & 1/2
\end{array}\right]^{\beta_j}$ if $\beta_j>0$ and
$\left[\begin{array}{cc}1/2 & -1/2 \\ -1/2 & 1/2+\epsilon
\end{array}\right]^{-\beta_j}$ if $\beta_j<0$, $j=1,\ldots,N$.

Thus, the limit of (\ref{epsilon}) for $\epsilon\to 0$ exists and
equals\begin{equation}\label{limit2} P_1Q_1P_2Q_2\cdots
P_NQ_N,\end{equation} where $P_j$ is $P=\left[\begin{array}{cc}1&0\\
0 & 0\end{array}\right]$ if $\alpha_j>0$ and $I-P$ if
$\alpha_j<0$, and $Q_j$ is $Q=\left[\begin{array}{cc}1/2 & 1/2\\
1/2 & 1/2\end{array}\right]$ if $\beta_j>0$ and $I-Q$ if
$\beta_j<0$.

A straightforward computation shows that \begin{gather*}
PQP=P(I-Q)P=\frac{1}{2} P,\
(I-P)Q(I-P)=(I-P)(I-Q)(I-P)=\frac{1}{2}(I-P),\\
\intertext{and}QPQ=Q(I-P)Q=\frac{1}{2} Q,
(I-Q)P(I-Q)=(I-Q)((I-P)(I-Q)=\frac{1}{2}(I-Q).\end{gather*}
Consequently (recall the condition imposed on the signs of
$\alpha$'s and $\beta$'s and the alternating nature of irreducible
sequences), the matrix (\ref{limit2}), up to a positive scalar
multiple $2^{N-m}$, coincides with \[ (PQ(I-P)(I-Q))^{m/2}\text{
if }m\text{ is even, and }(PQ(I-P)(I-Q))^{(m-1)/2}PQ \text{ if }m
\text{ is odd.}\]

It can be checked by induction that, for any $k\in\N$, \[
(PQ(I-P)(I-Q))^k=\frac{1}{4^k}\left[\begin{array}{cc} (-1)^k &
(-1)^{k-1}\\ 0 & 0\end{array}\right] .\] This implies that the
trace of (\ref{limit2}) in case of odd $\lfloor m/2\rfloor$ is
negative. But then, for sufficiently small $\epsilon>0$, the trace
of $A^{\alpha_1}B^{\beta_1}A^{\alpha_2}B^{\beta_2}\cdots
A^{\alpha_N}B^{\beta_N}$ also is negative. It remains to observe
that $\lfloor m/2\rfloor$ is odd if and only if $m\equiv 2$ or $3$
mod $4$.
\end{proof}

Theorems~\ref{th:2eig} and \ref{th:not2good} combined give a
complete description of all 2-good sequences of class 2. Observe
that in this case ``2-goodness'' does not depend on the magnitude
of the elements of the sequence but only on the sign pattern. At
the moment, we do not know whether this is true for sequences of
arbitrary length. We observe also that representation (\ref{proj})
shows that, if the sequence $\alpha_1,\ldots,\beta_N$ in
(\ref{word}) is 2-good, then $W(A,B)$ has positive eigenvalues for
matrices $A,B$ of {\em any} size, provided that each of them has
at most two distinct eigenvalues.
\providecommand{\bysame}{\leavevmode\hbox
to3em{\hrulefill}\thinspace}
\providecommand{\MR}{\relax\ifhmode\unskip\space\fi MR }
\providecommand{\MRhref}[2]{%
  \href{http://www.ams.org/mathscinet-getitem?mr=#1}{#2}
} \providecommand{\href}[2]{#2}

\end{document}